\documentclass[a4paper,10pt]{article}
\usepackage{amsmath,latexsym,amsfonts,amssymb,amsthm,mathrsfs,longtable,lscape,dsfont,fancybox,yfonts}
\usepackage{amsthm}
\usepackage{amssymb}
\usepackage{amscd}
\usepackage{amsfonts}
\usepackage{amsbsy}
\usepackage{graphicx}
\usepackage[all]{xy}
\usepackage{epstopdf}
\newtheorem {teo} {Theorem} [section]

\newtheorem {dfn} [teo] {Definition}
\date{}

\title{On the``Blue sky catastrophe" termination in the restricted four-body problem}
\author{Jaime Burgos--Garc\'ia \thanks{Departamento de Matem\'aticas
UAM--Iztapalapa. Av. San Rafael Atlixco 186, Col. Vicentina, C.P.
09340, M\'exico, D.F. e--mail: jbg84@xanum.uam.mx} \and Joaqu\'in
Delgado\thanks{Departamento de Matem\'aticas UAM--Iztapalapa.
e--mail: jdf@xanum.uam.mx}\and Submitted to Celestial Mechanics and Dynamical Astronomy}
\begin{document}

\maketitle

\begin{abstract}
The restricted three-body problem posses the property that some classes of doubly asymptotic orbits are limits members of families of periodic orbits, this phenomena has been known as the "Blue Sky Catastrophe" termination. A similar case occurs in the restricted four body problem for the collinear equilibrium point named $L_{2}$. We make an analytical and numerical study of the stable and unstable manifolds to verify that the hypothesis under which this phenomena occurs are satisfied.
\end{abstract}

\section{Introduction}
Few bodies problems have been studied for long time in celestial
mechanics,  either as simplified models of more complex planetary
systems or as benchmark models where new mathematical theories can
be tested. The three--body problem has been source of inspiration
and study in Celestial Mechanics since Newton and Euler, in particular the restricted three body problem (R3BP) has
demonstrated to be a good model of several systems in our solar
system such as the Sun-Jupiter-Asteroid system, and with less
accuracy the Sun-Earth-Moon system.  In analogy with the R3BP, in
this paper we study a restricted problem of four bodies consisting
of three primaries moving in circular orbits keeping an equilateral
triangle configuration and a massless particle moving under the
gravitational attraction of the primaries.
In 1933 in the Copenhagen observatory, Elis Str\"{o}mgren and co-workers performed a numerical investigation of the periodic orbits of the R3BP for the value $\mu=1/2$ of the mass parameter. In this investigation we can find a family of periodic orbits called \textit{class g}. This family originates form periodic orbits which Poincar\'e called \textit{premi\'ere sorte}, \cite{HenGen}, \cite{Sz}. From this simple beginning, the family develops many variations and the natural end of this family could not be stated by the investigators in the Copenhagen observatory. Str\"{o}mgren wrote in his conclusions that the termination of this family had to be an asymptotic periodic orbit spiralling into the equilibrium points $L_{4}$ and $L_{5}$. Several subsequent numerical explorations were necessary to confirm this conjecture. However, an analytical proof of this conjecture was given by J. Henrard \cite{Henr} until 1972 and later by Buffoni \cite{Buf}. In those papers, Henrard and Buffoni give sufficient conditions to prove the conjecture in the more general framework of analytical Hamiltonian systems with two degrees of freedom. These conditions amount to a transversality condition on the stable and unstable manifolds.
In papers like \cite{Burgos} \cite{PapaI} there exist numerical explorations of several families of periodic orbits of the restricted four--body problem for some values of the masses, in these works we can find that some families become asymptotic to a equilibrium point named $L_{2}$ as the value of the Jacobi constant at the equilibrium point is reached. So as in the R3BP, this numerical evidence requires to perform an analytical investigation to verify that the conditions for which the ``Blue sky catastrophe" termination ocurrs are satisfied.
This paper is organized as follows: In section 2 we state the equations of motion of the restricted four--body problem and we give some basic properties of the problem. In section 3 we analyse the linealization of the equations of motion at the equilibrium point $L_{2}$, in particular we will analyse the eigenvalues of the linear equations. In section 4 we give a brief introduction to the normal form theory and an application to the restricted four--body problem will be performed. In section 5 we analyse the normal form obtained and we will verify the conditions given in \cite{MeyMcSw} to prove a topologically transverse intersection between the stable and unstable manifolds.

\section{Equations of Motion}
\begin{figure}[!hbp]
\centering
\includegraphics[width=2.5in]{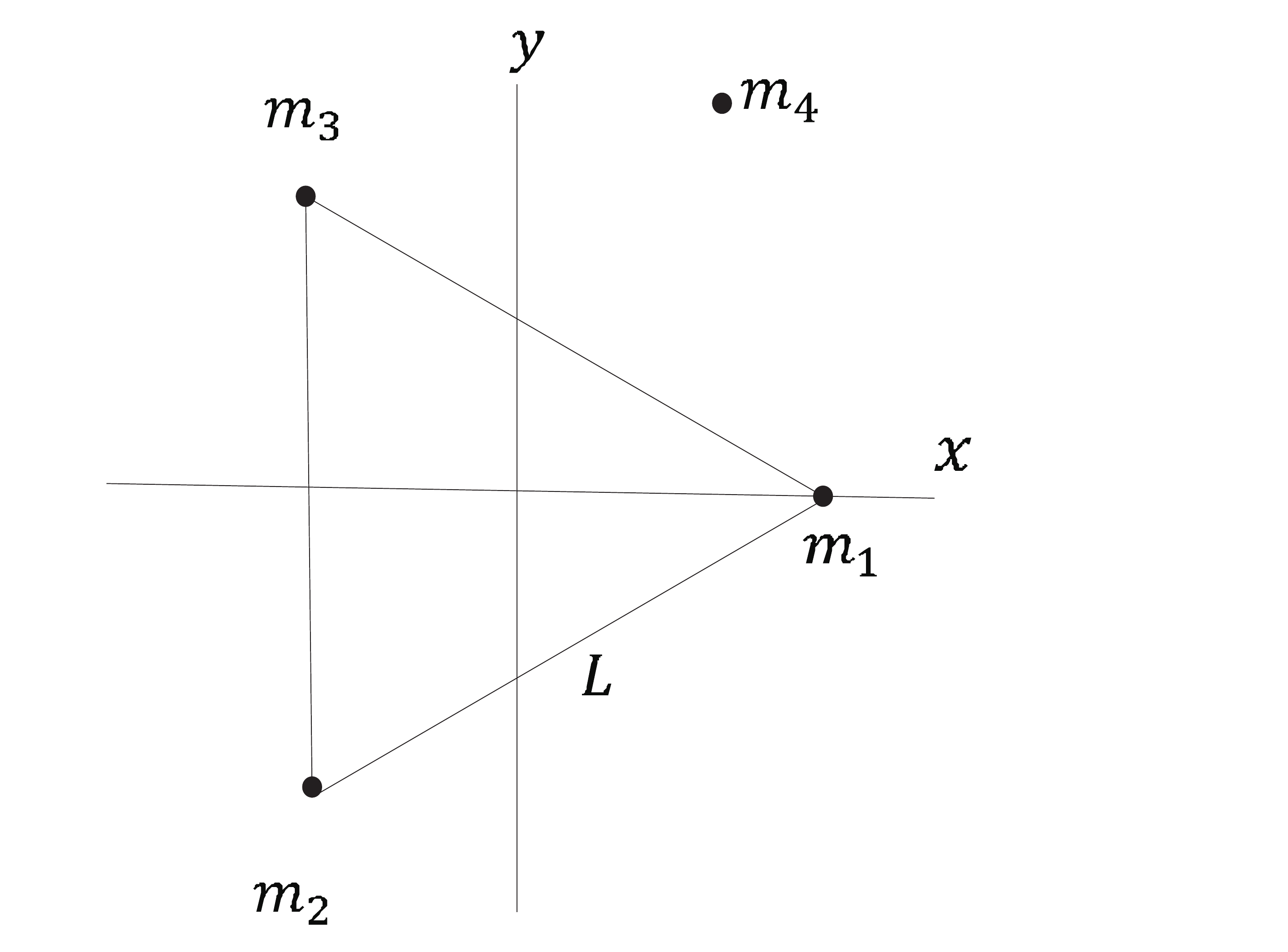}
\caption{The restricted four-body problem in a synodic system\label{triangle}}
\end{figure}

Consider three point masses, called $\textit{primaries}$, moving in
circular periodic orbits around their center of mass under their
mutual Newtonian gravitational attraction,  forming an equilateral
triangle configuration. A third massless particle moving in the same
plane is acted upon the attraction of the primaries. The equations
of motion of the massless particle referred to a synodic frame with
the same origin, where the primaries remain fixed, are:
\begin{eqnarray}
\bar{x}''-2n\bar{y}'-n^2\bar{x}&=&-k^2\sum_{i=1}^{3}m_{i}\frac{(\bar{x}-\bar{x_{i}})}{\rho_{i}^3}\nonumber\\
\bar{y}''+2n\bar{x}'-n^2\bar{y}&=&-k^2\sum_{i=1}^{3}m_{i}\frac{(\bar{y}-\bar{y_{i}})}{\rho_{i}^3} \label{sistemaconunidades}
\end{eqnarray}
where $k^2$ is the gravitational constant, $n$ is the mean motion,
$\rho_{i}^{2}=(\bar{x}-\bar{x}_{i})^2+(\bar{y}-\bar{y}_{i})^2$ is the distance of the massless particle to the primaries,
 $\bar{x}_{i}$, $\bar{y}_{i}$
are the vertices of equilateral triangle  formed by the primaries,
and ($'$)  denotes derivative with respect to time $t^{*}$. We
choose the orientation of the triangle of masses such that $m_1$
lies along the positive $x$--axis and $m_2$, $m_3$ are located
symmetrically with respect to the same axis, see
figure~\ref{triangle}.

The equations of motion can be recast in dimensionless form as follows:
Let $L$ denote the length of triangle formed by the primaries,  $x=\bar{x}/L$,
$y=\bar{y}/L$, $x_i=\bar{x}_i/L$, $y_i=\bar{y}_i/L$, for $i=1,2,3$;
$M=m_{1}+m_{2}+m_{3}$ the total mass, and  $t=nt^{*}$. Then the equations (\ref{sistemaconunidades}) become
\begin{eqnarray}
\ddot{x}-2\dot{y}-x &=&-\sum_{i=1}^{3}\mu_{i}\frac{(x-x_{i})}{r_{i}^3}\nonumber\\
\ddot{y}+2\dot{x}-y &=&-\sum_{i=1}^{3}\mu_{i}\frac{(y-y_{i})}{r_{i}^3} \label{sistemasinunidades}
\end{eqnarray}
where we have used Kepler's third law: $k^2M=n^2L^3$, and the dot
($\dot{}$) represents derivatives with respect to the dimensionless
time $t$ and $r_{i}^2=(x-x_{i})^2+(y-y_{i})^2$.

 The system
(\ref{sistemasinunidades}) will be defined if we know the vertices
of triangle for each value of the masses. In this paper we suppose
$\mu:=\mu_{3}=\mu_{2}$ then $\mu_{1}=1-2\mu$, it's not hard to prove
that the vertices of triangle are given as function of the mass
parameter $\mu$ by $x_{1}=\sqrt{3}\mu$, $y_{1}=0$,
$x_{2}=-\frac{\sqrt{3}(1-2\mu)}{2}$, $y_{2}=-\frac{1}{2}$,
$x_{3}=-\frac{\sqrt{3}(1-2\mu)}{2}$, $y_{3}=\frac{1}{2}$. The system
(\ref{sistemasinunidades}) can be written succinctly as

\begin{eqnarray}
\ddot{x}-2\dot{y}&=&\Omega_{x} \label{sistemastandar}\\
\ddot{y}+2\dot{x}&=&\Omega_{y}
\end{eqnarray}
where
\begin{displaymath}
\label{omega}\Omega(x,y,\mu):=\frac{1}{2}(x^2+y^2)+\sum_{i=1}^{3}\frac{\mu_{i}}{r_{i}}.
\end{displaymath}
is the effective potential function.

 There are three limiting cases:
\begin{enumerate}
\item If $\mu=0$, we obtain the rotating Kepler's problem, with $m_{1}=1$ at the origin of coordinates.
\item If $\mu=1/2$, we obtain the circular restricted three body problem, with two equal masses $m_{2}=m_{3}=1/2$.
\item If $\mu=1/3$, we obtain the symmetric case with three masses equal to $1/3$.
\end{enumerate}

It will be useful to write the system (\ref{sistemastandar}) using
complex notation. Let $z=x+\textit{i}y$, then
\begin{equation}
\label{sistemacomplejo}\ddot{z}+2\textit{i}\dot{z}=2\frac{\partial\Omega}{\partial\bar{z}}
\end{equation}
with
\begin{displaymath}
\Omega(z,\bar{z},\mu)=\frac{1}{2}\vert
z\vert^2+U(z,\bar{z},\mu)
\end{displaymath}
where the gravitational potential is
\begin{displaymath}
U(z,\bar{z},\mu)=\sum_{i=1}^{3}\frac{\mu_{i}}{\vert z-z_{i}\vert}
\end{displaymath}
and $r_{i}=\vert z-z_{i}\vert$, $i=1,2,3.$ are the distances to the
primaries. System (\ref{sistemacomplejo}) has the Jacobian first
integral
\begin{equation}
\label{integralprimera} 2\Omega(z,\bar{z},\mu)-\vert\dot{z}\vert^{2}=C
\end{equation}

If we define $P=p_{x}+\textit{i}p_{y}$, the conjugate momenta of
$z$, then system (\ref{sistemastandar}) can be recast as  a
Hamiltonian system  with Hamiltonian
\begin{eqnarray}
H&=&\frac{1}{2}\vert P\vert^2+Im(z\overline{P})-U(z,\bar{z},\mu)\nonumber\\
&=&\frac{1}{2}(p^{2}_{x}+p^{2}_{y})+(yp_{x}-xp_{y})-U(x,y,\mu).\label{hamiltoniano}
\end{eqnarray}
The relationship with the Jacobian integral is $H=-C/2$.
The phase space of (\ref{hamiltoniano}) is defined as
\begin{displaymath}
\Delta=\{(z,P)\in\mathbb{C}\times\mathbb{C}\vert z\ne z_{i},
i=1,2,3\},
\end{displaymath}
with collisions occurring at $z=z_{i}$, $i=1,2,3$.

\begin{figure}
  \centering
  \includegraphics[width=1.8in]{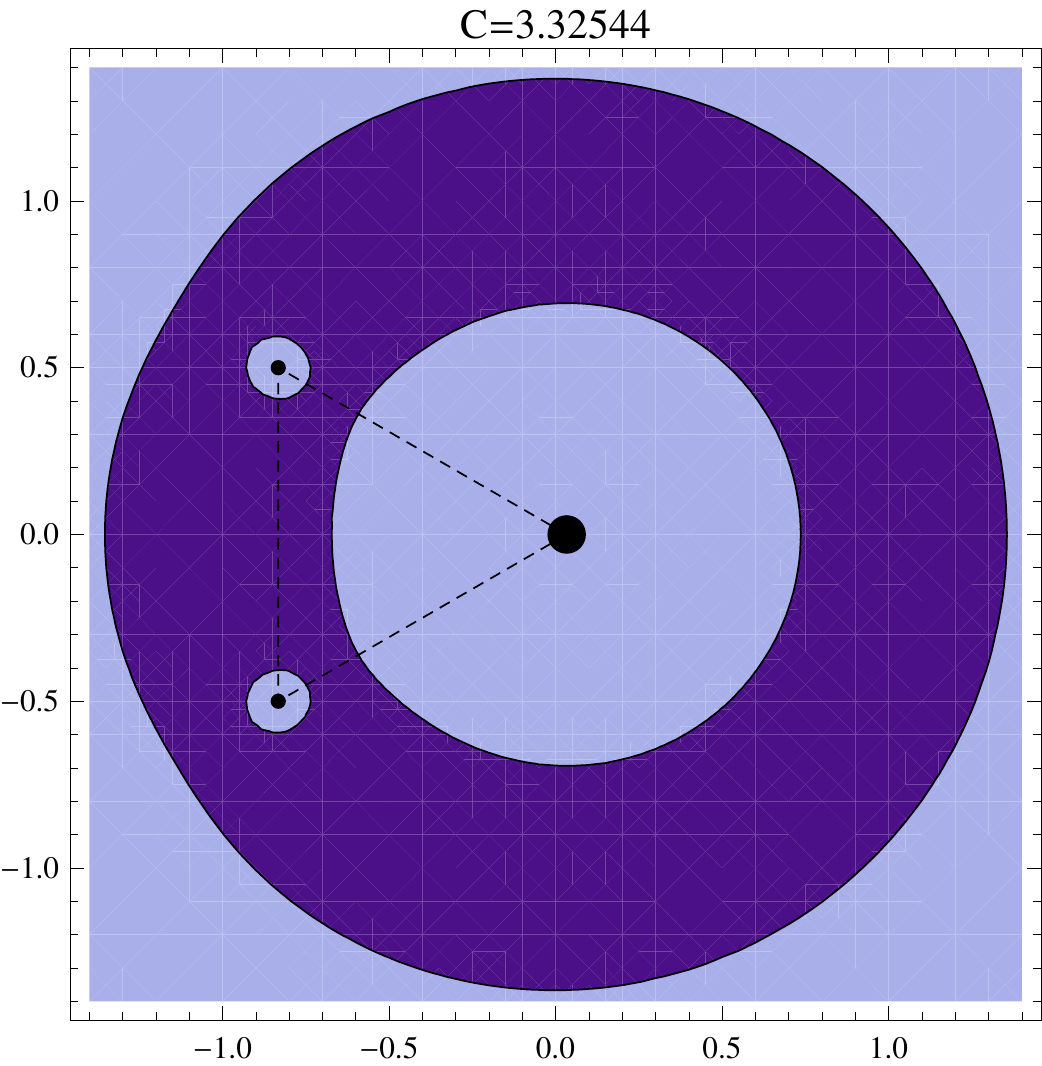}\quad
  \includegraphics[width=1.8in]{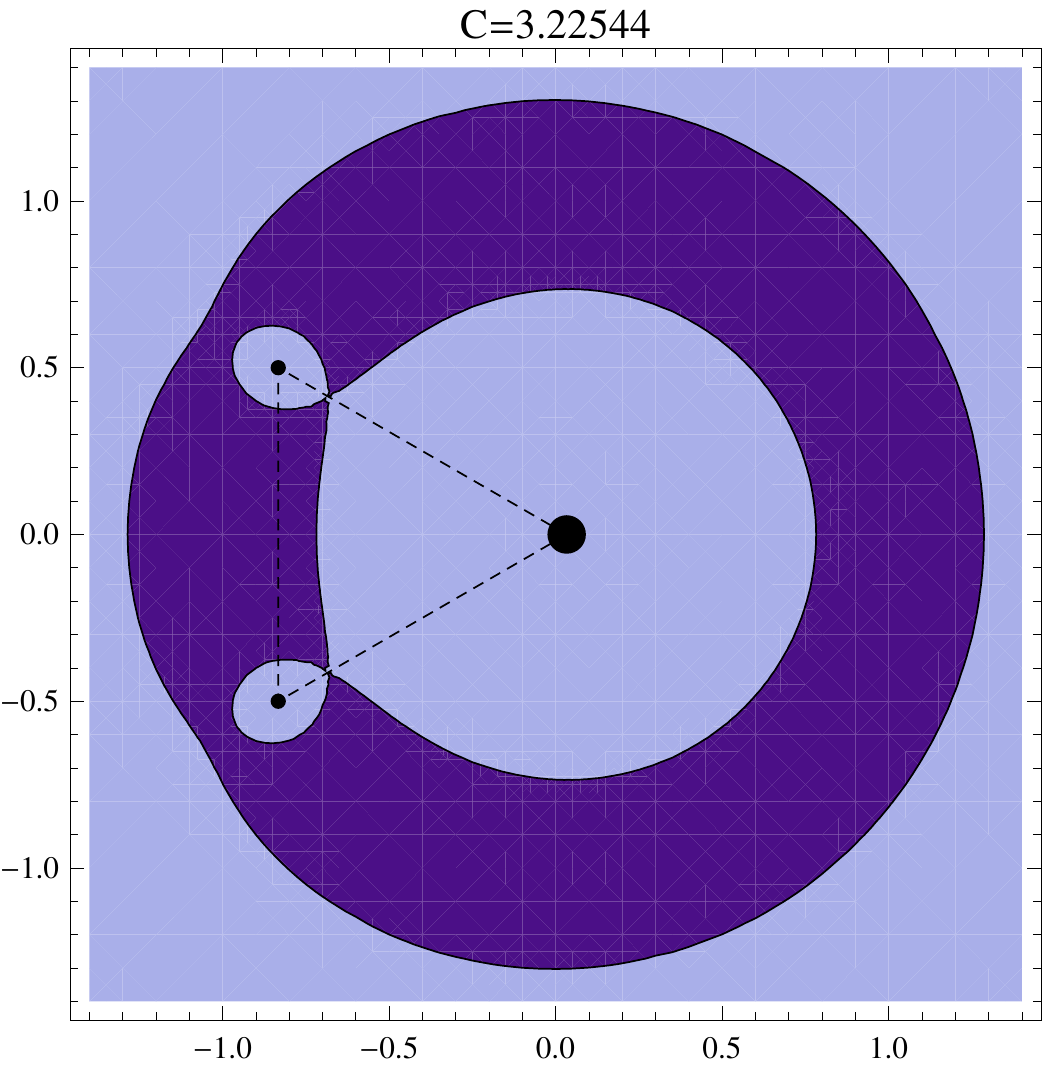}\\
  \includegraphics[width=1.8in]{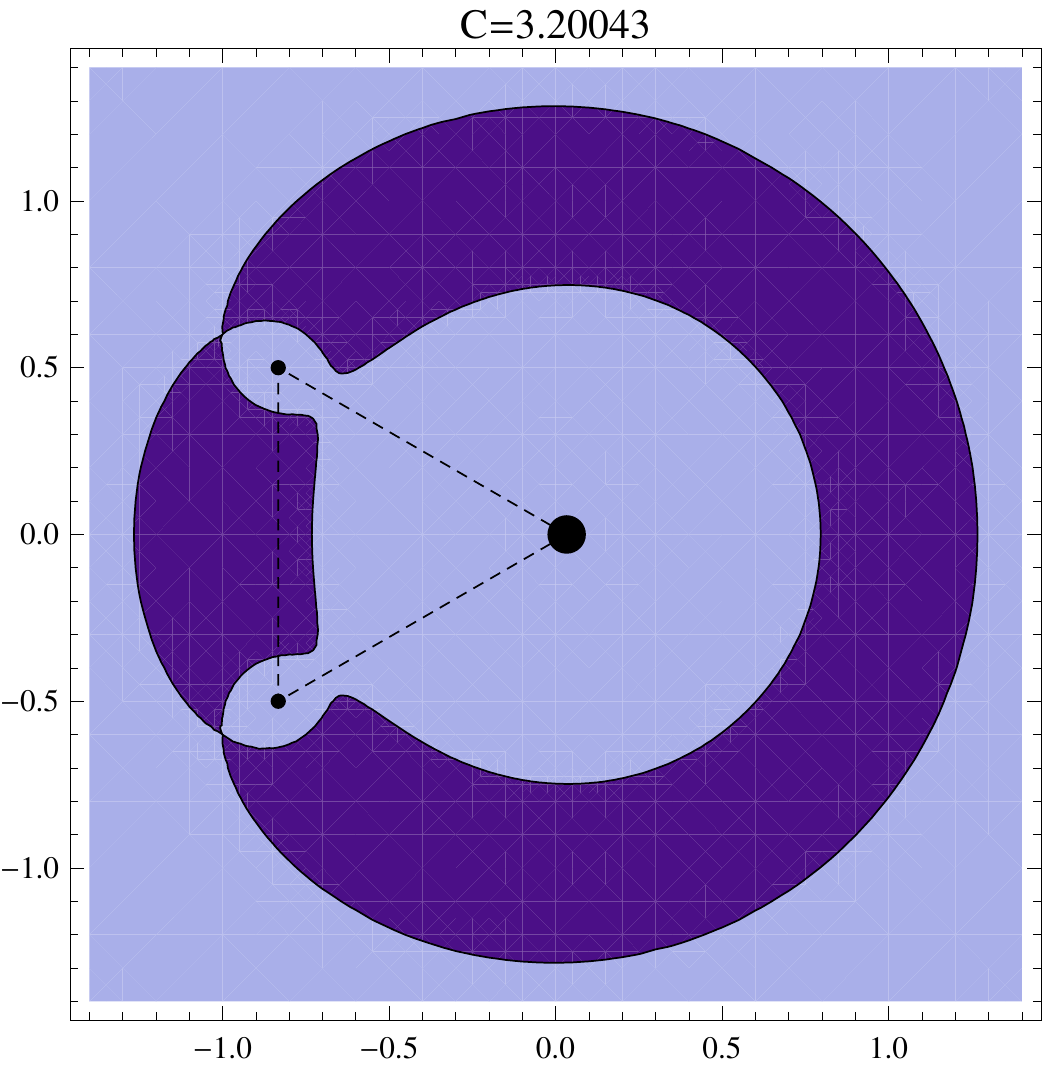}\quad
  \includegraphics[width=1.8in]{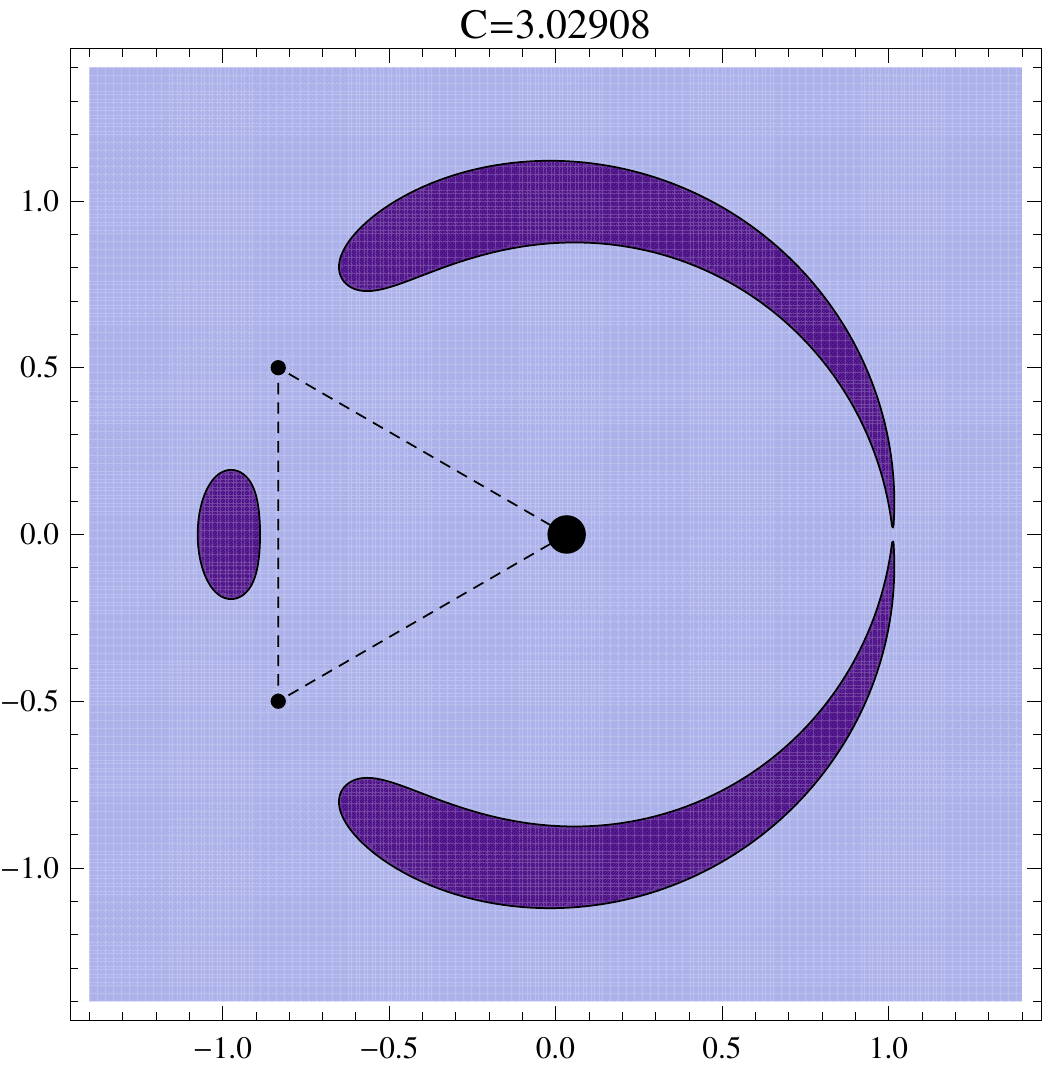}\\
  \includegraphics[width=1.8in]{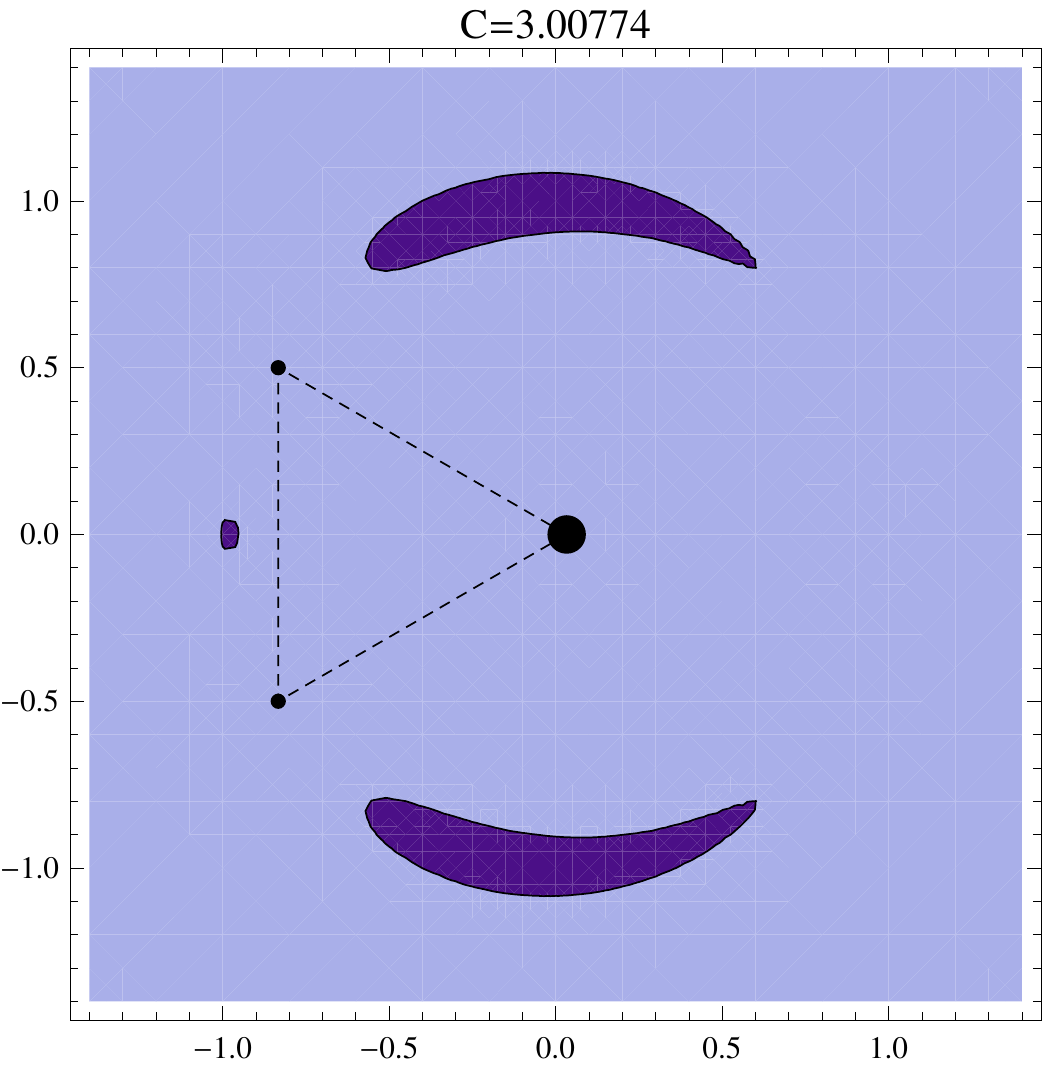}\quad
  \includegraphics[width=1.8in]{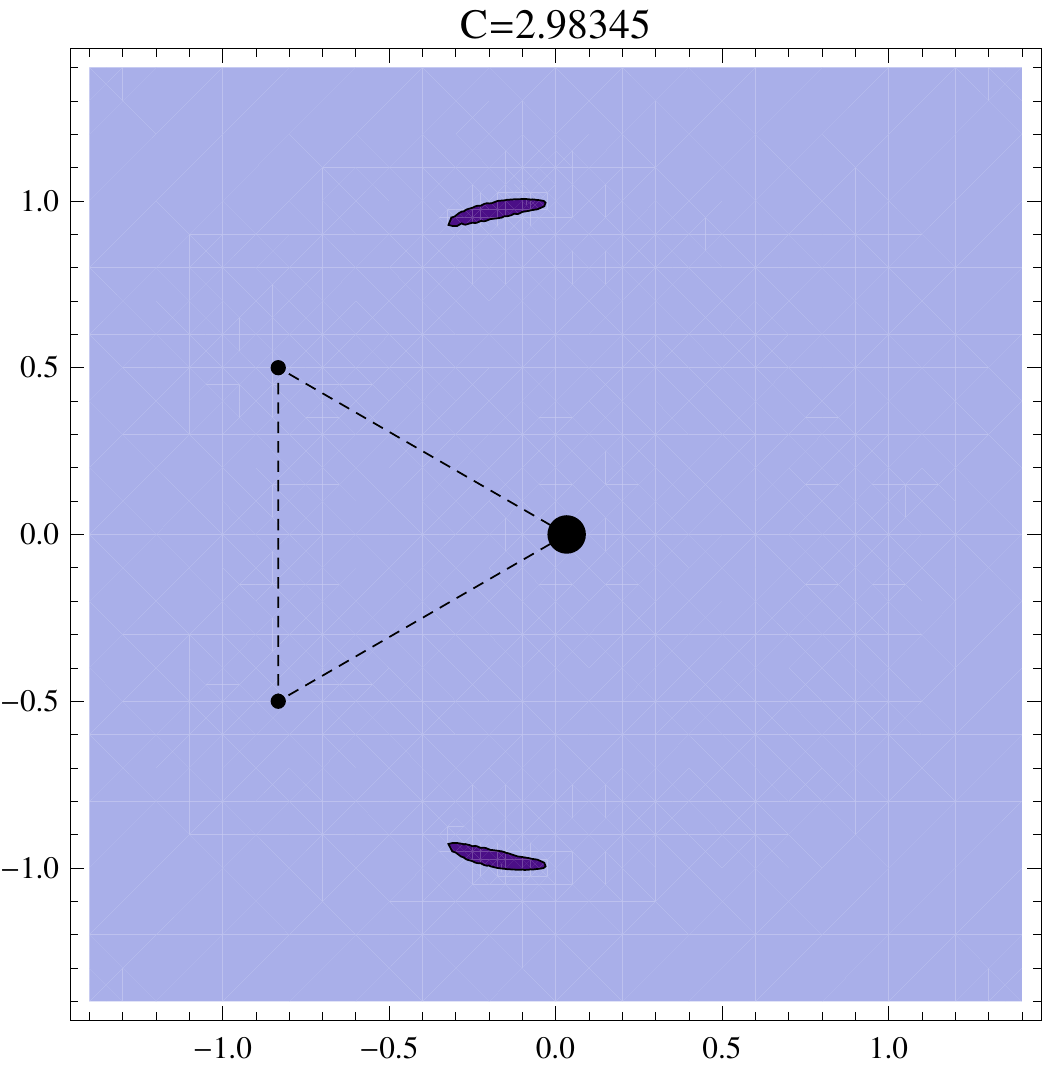}\\
      \caption{Hill's regions for a large value of the Jacobian constant (top-left). Hill's regions for critical values of the Jacobian constant. The Hill's regions of the last row correspond to a slightly larger value than the critical one for illustrations purposes. }\label{others}
\end{figure}

In the restricted three-body problem there exist five equilibrium
points for all values of the masses of the primaries but in this restricted
four-body problem (R4BP) the number of equilibrium
points depends on the particular values of the masses. Figure~\ref{others} shows the Hill's regions for a large value of the Jacobian constant $C$
consisting of a large exterior region around the primaries and  three components containing them.
As the critical Jacobian constant decreases, the evolution of the Hill's region is shown in Figure~\ref{others}. The smaller value of $C$ is just above the critical value where the Hill's region is the whole plain minus the positions of the primaries.

\begin{figure}[!h]
\begin{center}
\includegraphics[width=2.5in]{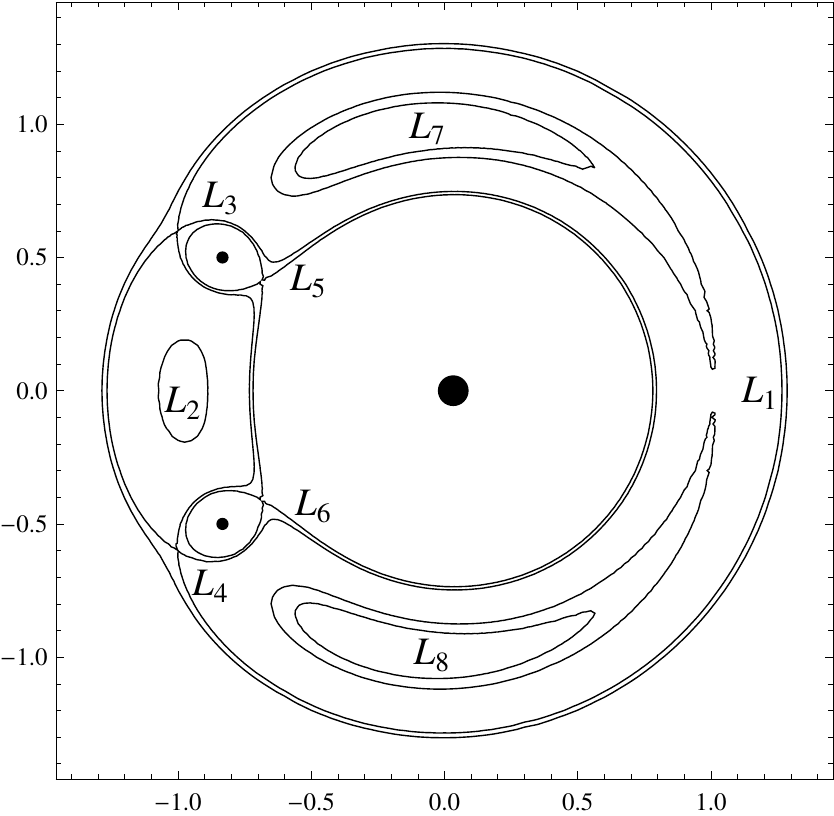}
\end{center}
\caption{The eight equilibrium points for $\mu=0.019$.\label{eqpoints}}
\end{figure}

A complete discussion of the equilibrium points and
bifurcations can be found in \cite{Del}, \cite{MeyerCC}, \cite{Lea},  \cite{PapaII}, \cite{Simo}.
For the particular value of the mass parameter presented in \cite{Burgos}, the R4BP has 2 collinear and 6 non-collinear equilibrium points.  We use the notation shown in Figure~\ref{eqpoints} for the eight critical points.

\section{Linearization around the collinear equilibrium point $L_{2}$}
In the former section it was shown that the number of equilibrium points in the restricted four-body problem depend on the value of the masses of the primaries. In papers like \cite{Lea} and \cite{Simo} we can see that the stability of the collinear equilibrium points also depends on the values of the masses, in this paper we will perform a thorough study of the stability of the equilibrium point $L_{2}$ (see Figure \ref{eqpoints}) and its consequences.

If we write $L_{2}=(\xi,0)$ for the coordinates of the equilibrium point then the Hamiltonian function centred at $L_2$ reads
\begin{displaymath}
\label{hamiltonianoreal} H=\frac{1}{2}(y^{2}_{1}+y^{2}_{2})+(x_{2}y_{1}-x_{1}y_{2})-x_{1}\xi-U(x_{1},x_{2},\mu)
\end{displaymath}
with
\begin{displaymath}
U(x_{1},x_{2},\mu)=\sum_{i=1}^{3}\frac{\mu_{i}}{r_{i}}
\end{displaymath}
where $r^{2}_{i}=(x_{1}+\xi-u_{1})^{2}+(x_{2}-v_{i})^{2}$, $i=1,2,3.$ and $u_{i}$, $v_{i}$ denote the positions of the primaries. Expanding through second-order terms, we obtain
\begin{equation}
\label{expandhamilt}
H=\frac{1}{2}(y^{2}_{1}+y^{2}_{2})+(x_{2}y_{1}-x_{1}y_{2})-\frac{1}{2}\left(U_{x_{1}x_{1}}x^{2}_{1}+U_{x_{2}x_{2}}x^{2}_{2}+......\right)
\end{equation}
There are no linear terms because the expansion is performed at a equilibrium point and the constant term has been omitted because it does not contribute to the system of differential equations. The term $U_{x_{1}x_{2}}$ is zero because the equilibrium point is collinear. If we ignore the higher order terms, the corresponding quadratic Hamiltonian gives rise the following Hamiltonian matrix
\begin{displaymath}
A=\left(\begin{array}{cccc}
0 & 1 & 1 & 0\\
-1 & 0 & 0 & 1\\
a & 0 & 0 & 1\\
0 & b & -1 & 0\\
\end{array}\right)
\end{displaymath}
Where $a=U_{x_{1}x_{1}}$ and $b=U_{x_{2}x_{2}}$, the characteristic equation of this matrix is
\begin{displaymath}
P(\lambda)=\lambda^{4}+(2-a-b)\lambda^{2}+(ab+a+b+1)
\end{displaymath}
using the relation $\Omega_{x_{1}x_{1}}=1+U_{x_{1}x_{1}}$, $\Omega_{x_{2}x_{2}}=1+U_{x_{2}x_{2}}$ we can write as in the R3BP:
\begin{equation}\label{charpoly}
P(\lambda)=\lambda^{4}+(4-\Omega_{x_{1}x_{1}}-\Omega_{x_{2}x_{2}})\lambda^{2}+\Omega_{x_{1}x_{1}}\Omega_{x_{2}x_{2}}.
\end{equation}
The analysis made in \cite{Lea} states that at $L_{2}$ there exist $\mu_{b}$ such that for $\mu<\mu_{b}$ the eigenvalues are $\pm\textit{i}\omega_{1}$, $\pm\textit{i}\omega_{2}$, for $\mu=\mu_{b}$ is $(\pm\textit{i}\omega)^{2}$ with multiplicity two, and for $\mu>\mu_{b}$ the eigenvalues are off the imaginary axis. The value $\mu_{b}$ can be determined numerically as a  zero of the discriminant of
\begin{equation}
\label{discriminant}
D=(4-\Omega_{x_{1}x_{1}}-\Omega_{x_{2}x_{2}})^{2}-4\Omega_{x_{1}x_{1}}\Omega_{x_{2}x_{2}}
\end{equation}
obtained form (\ref{charpoly}) if we make $\eta=\lambda^{2}$. This value is approximately $\mu_{b}=0.0027$, it must be clear that $\omega$ is in function of $a$ and $b$ and therefore of $\mu$. As we said in section 1, there exists numerical evidence that some families of periodic orbits of this problem undergo a  ``Blue Sky Catastrophe" (BSC) termination, therefore this suggest to verify the hypothesis given in \cite{Henr} and \cite{Buf} to give a rigorous demonstration of these phenomena.

\section{Normal form at the equilibrium point}
Perturbation theory is one of the few ways to study the behaviour of a real nonlinear system beyond its linear approximation, there exist a lot of papers dealing with this problem, for example in \cite{vander2}, \cite{wigg}, \cite{Poin} \cite{MeyerHDS} and \cite{Burgoy} it can be found a good introduction to the perturbation theory. However, we show briefly Deprit's algorithm in order to calculate a normal form at the equilibrium point $L_{2}$.

Let $H(\epsilon,x)$ be a Hamiltonian, $X(\epsilon,y)$ be a change of variables generated by a function $W(\epsilon,x)$ and $G(\epsilon,y)=H(\epsilon,X(\epsilon,y))$. Suppose that $H$, $W$, $G$ all have series expansions in the small parameter $\epsilon$. The algorithm of the method of Lie transforms is a recursive set of formulas that relate the terms in these various series expansions i,e., let
\begin{equation}
\label{originalseries}
H(\epsilon,x)=H_{0}^{0}(x)+\sum_{n=1}^{\infty}\left(\frac{\epsilon^{n}}{n!}\right) H_{n}^{0}(x)
\end{equation}
\begin{equation}
\label{transforseries}
G(\epsilon,y)=H_{0}^{0}(y)+\sum_{n=1}^{\infty}\left(\frac{\epsilon^{n}}{n!}\right) H_{0}^{n}(x)
\end{equation}
\begin{equation}
\label{generatingfun}
W(\epsilon,x)=\sum_{n=0}^{\infty}\left(\frac{\epsilon^{n}}{n!}\right) W_{n+1}(x)
\end{equation}
the functions $H_{j}^{i}$, $i=1,2,..$, $j=0,1,2...$ satisfy the recursive identities
\begin{displaymath}
\label{recursiveformula}
H_{j}^{i}=H_{j+1}^{i-1}+\sum_{k=0}^{j}{j \choose k}\{H_{j-k}^{i-1},W_{k+1}\}
\end{displaymath}
For example, to compute the series $G(\epsilon,y)$ through terms of order $\epsilon^{2}$ we need to solve the first \textit{homological equation}
\begin{equation}
\label{homological1}
H_{0}^{1}=H_{1}^{0}+\{H_{0}^{0},W_{1}\}
\end{equation}
which gives the term of order $\epsilon$, then we must  compute
\begin{displaymath}
H_{1}^{1}=H_{2}^{0}+\{H_{1}^{0},W_{1}\}+\{H_{0}^{0},W_{2}\}
\end{displaymath}
and  solve a second homological equation
\begin{equation}
\label{homological2}
H_{0}^{2}=H_{1}^{1}+\{H_{0}^{1},W_{1}\}
\end{equation}
which gives the term of order $\epsilon^{2}$ in this way we  obtain
\begin{equation}
\label{transforseries2}
G(\epsilon,y)=H_{0}^{0}(y)+\epsilon H_{0}^{1}(y)+\frac{\epsilon^{2}}{2}H_{0}^{2}(y)+....
\end{equation}
We can transform the Hamiltonian (\ref{expandhamilt}) in the Hamiltonian (\ref{originalseries}) if we make the symplectic scaling $x=\epsilon(x_{1},x_{2},y_{1},y_{2})\in\mathbb{R}^{4}$ with multiplier $\epsilon^{-2}$ in the Hamiltonian (\ref{expandhamilt}) and collect the homogeneous terms.

\subsection{Linear normal form}
If we ignore the higher order terms in (\ref{expandhamilt}) we get the following quadratic Hamiltonian
\begin{equation}
\label{quadhamil}
H=\frac{1}{2}(y^{2}_{1}+y^{2}_{2})+(x_{2}y_{1}-x_{1}y_{2})-\frac{1}{2}\left(U_{x_{1}x_{1}}x^{2}_{1}+U_{x_{2}x_{2}}x^{2}_{2}\right),
\end{equation}
this Hamiltonian can be written as a quadratic form $H=\frac{1}{2}x^{T}Sx$, with $S=-J_{4}A$, $J_{4}$ denotes the standard four dimensional symplectic matrix and $A$ is the matrix defined in the previous section. We seek a symplectic change of coordinates $x=Pz$ such that the Hamiltonian (\ref{quadhamil}) can be written as $H=\frac{1}{2}z^{T}S^{*}z$ where $S^{*}=-J_{4}B$ and $B=P^{-1}AP$ is the normal form of $A$
\begin{equation}
\label{normalmatrix}
B=\left(\begin{array}{cccc}
0 & -\omega & 0 & 0\\
\omega & 0 & 0 & 0\\
\epsilon & 0 & 0 & -\omega\\
0 & \epsilon & \omega & 0\\
\end{array}\right)
\end{equation}
where $\epsilon=\pm1$, this value depends on the basis. In papers like \cite{kock} and \cite{Burgoy} there are detailed discussions on how to find the matrix $P$, in particular we follow the algorithm shown in \cite{Burgoy} to achieve this. The symplectic basis required depends on the decomposition $A=\Sigma+N$ where $\Sigma$ is a real semisimple symplectic matrix and $N$ is nilpotent matrix given by $N=A-\Sigma$. In terms of the quantities $a$, $b$ and $\omega$ as in the R3BP these matrices looks like $$\Sigma=\frac{1}{2\omega^{2}}$$
\begin{equation}
\label{semi}
\left(\begin{array}{cccc}
0 & 3\omega^{2}+2b+a-1 & 3\omega^{2}+a-3 & 0\\
-(3\omega^{2}+2a+b-1) & 0 & 0 & 3\omega^{2}+b-3\\
a^{2}-b+a(3\omega^{2}-2) & 0 & 0 & 3\omega^{2}+2a+b-1\\
0 & -a+b(3\omega^{2}+b-2) & -(3\omega^{2}+2b+a-1) & 0\\
\end{array}\right)
\end{equation}
and $$N=\frac{1}{2\omega^{2}}$$
\begin{equation}
\label{nil}
\left(\begin{array}{cccc}
0 & -(\omega^{2}+2b+a-1) & -(\omega^{2}+a-3) & 0\\
\omega^{2}+2a+b-1 & 0 & 0 & -(\omega^{2}+b-3)\\
-(a^{2}-b+a(\omega^{2}-2)) & 0 & 0 & -(\omega^{2}+2a+b-1)\\
0 & a-b(\omega^{2}+b-2) & \omega^{2}+2b+a-1 & 0\\
\end{array}\right)
\end{equation}
We must observe that the structures of $\Sigma$ and $N$ are equal to those of the R3BP in the non semisimple case, more precisely, if we substitute the numerical values of $a$ and $b$ we can see that $N\neq0$ but $N^{2}=0$ and it follows that the nilpotent index equals to 1; therefore we can build the desired basis as it was done in the R3BP. First consider the symplectic product on $\mathbb{R}^{4}$, $\langle x,y\rangle=x^{T}J_{4}y$, second, we need a initial vector to build the basis and to determine the value of $\epsilon$. As in the R3BP a natural choose is the vector $z_{0}=\frac{e_{1}}{\sqrt{\vert\langle e_{1},Ne_{1}\rangle\vert}}$ with $e_{1}$ the standard unitary vector in $\mathbb{R}^{4}$, it is not hard to see that $\langle e_{1},Ne_{1}\rangle=N_{31}$ and therefore
\begin{displaymath}
\epsilon=\langle z_{0},Nz_{0}\rangle=\frac{\langle e_{1},Ne_{1}\rangle}{\vert N_{31}\vert}=\frac{N_{31}}{\vert N_{31}\vert}
\end{displaymath}
the value of $N_{31}$ is approximately $-1.82$ and therefore $\epsilon=-1$ as in the R3BP. The desired basis has the form : $z_{1}=z_{0}+\frac{\epsilon}{2\omega^{2}}\langle z_{0},\Sigma z_{0}\rangle N\Sigma z_{0}$, $z_{2}=\frac{1}{\omega}\Sigma z_{1}$, $z_{3}=\epsilon Nz_{1}$, $z_{4}=\frac{\epsilon}{\omega}\Sigma Nz_{1}$. Therefore the desired matrix is  $P=col(z_{1},z_{2},z_{3},z_{4})$, this matrix has the form
\begin{equation}
\label{changebasismatrix}
P=\left(\begin{array}{cccc}
p_{11} & 0 & 0 & p_{14}\\
0 & p_{22} & p_{23} & 0\\
0 & p_{32} & p_{33} & 0\\
p_{41} & 0 & 0 & p_{44}\\
\end{array}\right)
\end{equation}
of course the entries of this matrix depend of the entries of $\Sigma$ and $N$ (therefore of $a$, $b$). The Hamiltonian (\ref{originalseries}) under the change of coordinates $z=Px$ becomes
\begin{equation}
\label{z=pzseries}
H(\epsilon,z)=H_{0}^{0}(z)+\sum_{n=1}^{\infty}\left(\frac{\epsilon^{n}}{n!}\right) H_{n}^{0}(z)
\end{equation}
$H_{0}^{0}(z)$ is the linear Hamiltonian associated with the normal form $B$.

\subsection{The higher order terms}
In this section we deal with the higher order terms in the expansion of the Hamiltonian (\ref{originalseries}), recalling the expansion given in (\ref{expandhamilt}) we see that the terms $H_{n}^{0}(x)$ contain the terms of order $n+2$ of the Taylor series of the potential $U(x_{1},x_{2},\mu)$, therefore such terms have the form $H_{1}^{0}(x)=a_{3}x_{1}^{3}+b_{3}x_{1}^{2}x_{2}+c_{3}x_{1}x_{2}^{2}+d_{3}x_{2}^{3}$, $H_{2}^{0}(x)=a_{4}x_{1}^{4}+b_{4}x_{1}^{3}x_{2}+c_{4}x_{1}^{2}x_{2}^{2}+d_{4}x_{1}x_{2}^{3}+e_{4}x_{2}^{4}$, ...etc. Here the coefficients of the these homogeneous polynomials correspond to the respective derivatives of order $n+2$ in the expansion of $U(x_{1},x_{2},\mu)$ evaluated at the origin. Because the equilibrium point is collinear we have that $\frac{\partial U}{\partial x_{2}}(x_{1},0,\mu)=0$ and therefore $\frac{\partial^{n}\partial}{\partial x_{1}^{n}\partial x_{2}}U(x_{1},0,\mu)=0$ for all $n$. A straightforward calculation shows
\begin{displaymath}
\frac{\partial^{3}U}{\partial x_{2}^{3}}=\mu\left(\frac{9(x_{2}-1/2)}{r_{3}^{5}}-\frac{15(x_{2}-1/2)^{3}}{r_{3}^{7}}\right)+\mu\left(\frac{9(x_{2}+1/2)}{r_{2}^{5}}-\frac{15(x_{2}+1/2)^{3}}{r_{2}^{7}}\right)
\end{displaymath}
\begin{displaymath}
+(1-2\mu)\mu\left(\frac{9x_{2}}{r_{1}^{5}}-\frac{15x_{2}^{3}}{r_{1}^{7}}\right)
\end{displaymath}
if we evaluate in $x_{2}=0$ we see that $r_{2}=r_{3}$, then $\frac{\partial^{3}U}{\partial x_{2}^{3}}=0$ and this implies $\frac{\partial^{4}\partial }{\partial x_{1}\partial x_{2}^{3}}U=0$. Therefore we have $H_{1}^{0}(x)=a_{3}x_{1}^{3}+c_{3}x_{1}x_{2}^{2}$ and $H_{2}^{0}(x)=a_{4}x_{1}^{4}+c_{4}x_{1}^{2}x_{2}^{2}+e_{4}x_{2}^{4}$. In the rest of this section we will not show all of the calculations. The next step to calculate the desired normal form is to transform $H_{1}^{0}$ and $H_{2}^{0}$ under the change of coordinates $x=Pz$, with help of a algebraic manipulator this is not problem, the interesting task is how to calculate the homological equations (\ref{homological1}) and (\ref{homological2}). There exist several ways on how to calculate these equations, in \cite{vander}, \cite{Soko} and \cite{Deprit} we can find a detailed discussions on this problem. Several authors use complex variables to perform the calculations of the normal form, however, it is worth noting that sometimes the normal forms under complex symplectic transformations do not determine the real normal forms \cite{Burgoy}, therefore we use only \textit{real} symplectic transformations to do the calculations. In \cite{Palayan} there is a elegant method to calculate the homological equations in the case non-semisimple of a linearized vector field and its application to the R3BP at the equilibrium point $L_{4}$, we are going to explain it briefly.\\ \\
Let us consider the change to symplectic polar coordinates $(z_{1},z_{2},z_{3},z_{4})\rightarrow(r,\theta,R,\Theta)$ given by
\begin{displaymath}
z_{1}=r\cos\theta
\end{displaymath}
\begin{displaymath}
z_{2}=r\sin\theta
\end{displaymath}
\begin{displaymath}
z_{3}=R\cos\theta-\frac{\Theta}{r}\sin\theta
\end{displaymath}
\begin{displaymath}
z_{4}=R\sin\theta+\frac{\Theta}{r}\cos\theta
\end{displaymath}
in these coordinates the term $H_{0}^{0}$ becomes $H_{0}^{0}=\omega\Theta+\frac{1}{2}r^{2}$, by scaling time we can assume that $\omega=1$; therefore
\begin{equation}
\label{H0polar}
H_{0}^{0}=\Theta+\frac{1}{2}r^{2}
\end{equation}
of course we must transform again the terms $H_{1}^{0}$, $H_{2}^{0}$ under this change of coordinates. Let $\mathfrak{L}_{H_{0}^{0}}(W_{i})$ the linear operator defined by $\mathfrak{L}_{H_{0}^{0}}(W_{i})=\{W_{i},H_{0}^{0}\}$, $i\in\mathbb{N}$, explicitly $\mathfrak{L}_{H_{0}^{0}}(W_{i})=\{W_{i},\Theta\}+\frac{1}{2}\{W_{i},r^{2}\}=\frac{\partial W_{i}}{\partial\theta}-r\frac{\partial W_{i}}{\partial R}:=\mathfrak{L}_{S}+\mathfrak{L}_{N}$ with $\mathfrak{L}_{S}=\frac{\partial}{\partial\theta}$ and $\mathfrak{L}_{N}=-r\frac{\partial}{\partial R}$. The homological equation in every step looks like
\begin{displaymath}
\mathfrak{L}_{H_{0}^{0}}(W_{i})+H_{0}^{i}=\widetilde{H_{i}}
\end{displaymath}
where $\widetilde{H_{i}}$ groups all the terms of the previous steps $H_{j}^{i}$. Deprit \cite{Deprit} proved that, after splitting $\widetilde{H_{i}}=\widetilde{H_{i}}^{*}(r,-,R,\Theta)+\widetilde{H_{i}^{'}}(r,\theta,R,\Theta)$, one can choose $H_{0}^{i}=\widetilde{H_{i}}^{*}$ and to solve
\begin{equation}
\label{homologpolar}
\mathfrak{L}_{H_{0}^{0}}(W_{i})=\widetilde{H_{i}^{'}}
\end{equation}
with $\widetilde{H_{i}^{'}}=\widetilde{H_{i}}-\widetilde{H_{i}}^{*}$ in terms of operators we can write the equation (\ref{homologpolar}) as
\begin{displaymath}
\left(id+\mathfrak{L}_{S}^{-1}\mathfrak{L}_{N}\right)(W_{i})=\mathfrak{L}_{S}^{-1}(\widetilde{H_{i}^{'}})
\end{displaymath}
in general the operator $\left(id+\mathfrak{L}_{S}^{-1}\mathfrak{L}_{N}\right)(W_{i})$ is not invertible, but in our particular case we can calculate $W_{i}$ because all the terms of the right side of (\ref{homologpolar}) are periodic in $\theta$ and $\mathfrak{L}_{N}$ is nilpotent. Because of the linearity of $\mathfrak{L}_{H_{0}^{0}}$ we have
\begin{displaymath}
\left(id+\mathfrak{L}_{S}^{-1}\mathfrak{L}_{N}\right)^{-1}=id-\left(\mathfrak{L}_{S}^{-1}\mathfrak{L}_{N}\right)+\left(\mathfrak{L}_{S}^{-1}\mathfrak{L}_{N}\right)^{2}-\cdots
\end{displaymath}
this series is a finite sum because in general $\mathfrak{L}_{N}^{m}(p)=0$ if $p$ is a homogeneous polynomial of degree $m-1$. Therefore $W_{i}$ is taken as
\begin{displaymath}
W_{i}=\left(id+\mathfrak{L}_{S}^{-1}\mathfrak{L}_{N}\right)^{-1}\mathfrak{L}_{S}^{-1}(\widetilde{H_{i}^{'}})
\end{displaymath}
observe that we obtain a function totally periodic in $\theta$. For this reason, we can repeat the process at any order. In particular after two steps we obtain
\begin{displaymath}
G(\epsilon,y)=H_{0}^{0}(y)+\epsilon H_{0}^{1}(y)+\frac{\epsilon^{2}}{2}H_{0}^{2}(y)
\end{displaymath}
The generating function can be known totally in every step without additional assumptions and it is entirely polynomial. Now we want to apply the method explained above to calculate the Hamiltonian (\ref{transforseries2}) in the case of the restricted four--body problem. Let $H_{0}^{0}$ be the Hamiltonian (\ref{H0polar}), as above, the first order homological equation to solve is
\begin{displaymath}
\mathfrak{L}_{H_{0}^{0}}(W_{1})+H_{0}^{1}=H_{1}^{0}
\end{displaymath}
If $H_{1}^{0}=H_{1}^{0^{*}}(r,-,R,\Theta)+H_{1}^{0^{'}}(r,\theta,R,\Theta)$
and we choose $H_{0}^{1}=H_{1}^{0^{*}}$ then the homological equation to solve is
\begin{equation}
\label{homologicalfirstorder}
\mathfrak{L}_{H_{0}^{0}}(W_{1})=H_{1}^{0^{'}}
\end{equation}
with $H_{1}^{0^{'}}=H_{1}^{0}-H_{1}^{0^{*}}$. In polar coordinates the term $H_{1}^{0}$ looks like $$H_{1}^{0}=a_{3}\left(p_{11}r\cos\theta+p_{14}\left(R\sin\theta+\frac{\Theta}{r}\cos\theta\right)\right)^{3}+$$
$$c_{3}\left(p_{22}r\sin\theta+p_{23}\left(R\cos\theta-\frac{\Theta}{r}\sin\theta\right)\right)^{2}\left(p_{11}r\cos\theta+p_{14}\left(R\sin\theta+\frac{\Theta}{r}\cos\theta\right)\right)$$
observe that the terms with $\theta$ appear in all the terms of $H_{1}^{0}$ therefore $H_{1}^{0^{*}}=0$ and $H_{1}^{0}=H_{1}^{0^{'}}$. If we write $I=\int H_{1}^{0}d\theta$ we can see that $\mathfrak{L}_{N}^{4}(I)=0$, so the first order generating function is given by
\begin{equation}
W_{1}=I+r\left(\int\frac{\partial}{\partial R}d\theta\right)I+r^{2}\left(\int\frac{\partial}{\partial R}d\theta\right)^{2}I+r^{3}\left(\int\frac{\partial}{\partial R}d\theta\right)^{3}I
\end{equation}
where the exponents stand by the composition of operators. For the second order term we have to calculate $H_{1}^{1}=H_{2}^{0}+\{H_{1}^{0},W_{1}\}+\{H_{0}^{0},W_{2}\}$ and $H_{0}^{2}=H_{1}^{1}+\{H_{0}^{1},W_{1}\}$ but $H_{0}^{1}=0$ then $H_{0}^{2}=H_{1}^{1}=\widetilde{H_{2}}+\{H_{0}^{0},W_{2}\}$ with $\widetilde{H_{2}}=H_{2}^{0}+\{H_{1}^{0},W_{1}\}$, therefore the second order homological equation is $$\mathfrak{L}_{H_{0}^{0}}(W_{2})=\widetilde{H_{2}^{'}}$$ with $H_{0}^{2}=\widetilde{H_{2}^{*}}$ and $\widetilde{H_{2}}=\widetilde{H_{2}^{*}}+\widetilde{H_{2}^{'}}$, the second order homological equation is solved as before. Note that along this section we have performed all the calculations without using the numerical values of the coefficients of the matrix $P$ and the polynomials $H_{1}^{0}$ and $H_{2}^{0}$. The expressions of $W_{1}$ and $H_{0}^{2}$ obtained with help of the algebraic manipulator are too long and we won't show them. Now if we substitute the required numerical coefficients of the matrix $P$: $p_{11}=-03928$, $p_{14}=-0.7631$, $p_{22}=-0.9680$, $p_{23}=-1.8807$, $p_{32}=2.005$, $p_{33}=1.3490$, $p_{41}=0.8134$, $p_{44}=0.5474$ and of the polynomials ones $a_{3}=-0.962$, $c_{3}=1.370$, $a_{4}=-1.007$, $c_{4}=3.150$, $e_{4}=-0.4686$ we obtain
\begin{equation}
\label{secondorderterm}
H_{0}^{2}=h_{1}r^{4}-h_{2}R^{4}+h_{3}r^{2}R^{2}-h_{4}r^{2}\Theta +h_{5}R^{2}\Theta -h_{6}\Theta^{2}-h_{7}\frac{R^{2}\Theta^{2}}{r^{2}}+
\end{equation}
$$h_{8}\frac{\Theta^{3}}{r^{2}}-h_{9}\frac{\Theta^{4}}{r^{4}}$$ with $h_{1}=2.19104$, $h_{2}=1.41252$, $h_{3}=16.2535$, $h_{4}=8.35177$, $h_{5}=4.24874$, $h_{6}=0.00392$, $h_{7}=2.82504$, $h_{8}=4.24874$, $h_{9}=1.41252$. Therefore the Hamiltonian (\ref{transforseries2}) after two steps looks like $$G=H_{0}^{0}+\frac{\epsilon^{2}}{2}H_{0}^{2}$$
In the non-semisimple case the normal form is not unique but it depends on how we perform its construction, in \cite{vander} it can be found good discussion on this topic, however, a straightforward calculation shows that the superior order term (\ref{secondorderterm}) contains the terms of $H_{0}^{2}=(p_{1}^{2}+p_{2}^{2})[A(p_{1}^{2}+p_{2}^{2})+B(q_{1}p_{2}-q_{2}p_{1})+C(q_{1}^{2}+q_{2}^{2})]$ and $H_{0}^
{2}=c(q_{1}p_{2}-q_{2}p_{1})^2+d(q_{1}p_{2}-q_{2}p_{1})(p_{1}^{2}+p_{2}^{2})+e(p_{1}^{2}+p_{2}^{2})^{2})$ with $A$, $B$, $C$, $c$, $d$, $e$ constants, found in \cite{Soko} and \cite{MeyerHDS} in polar coordinates. This is a direct consequence of our construction of the normal form where we are not putting additional conditions on the generating function $W$ and using only real coordinates.

\section{Analysis of the normal form}
We are interested in making a local analysis of the truncated system and of the stable and unstable manifolds. First, we will outline the construction of a versal deformation (\cite{vander}, \cite{Arnold}) of the normal form obtained in the previous section. It will useful recalling the definition of a versal deformation
\begin{dfn} \label{def} A deformation $A(\lambda)$ of $A_{0}\in M(\mathbb{C}_{n})$ is called $versal$ if any deformation $B(\mu)$ of $A_{0}$ is equivalent to a deformation induced from $A$, i.e.,
$$B(\mu)=C(\mu)A(\phi(\mu))C^{-1}(\mu)$$ for some change of parameters $\phi:\Sigma\rightarrow\Lambda$.
\end{dfn}
This kind of deformation is useful to study the truncated system through resonance. We already said we do not have a explicit relation of the partial derivatives $a$ and $b$ with the mass parameter $\mu$ so we cannot give a explicit relation between the parameters of the versal deformation with the mass parameter $\mu$. Therefore we just give a general relation of the parameters of the versal deformation with the quantities $a$ and $b$. Consider the matrix $B$ related to the linear vector field $H_{0}^{0}$, if we want to construct a versal deformation of the matrix $B$ we must take the matrix family $B_{\nu}:=B+D_{\nu}$ where the matrix $D_{\nu}$ lies on the orthogonal complement of the orbit of $B$ \cite{Wiggins}, this is equivalent to ask that the matrix $D_{\nu}$ satisfies $[D_{\nu}^{*},B]=0$ where $D_{\nu}^{*}$ denotes the adjoin of $D_{\nu}$. The matrix $D_{\nu}=\nu_{1}e_{1}+v_{2}e_{2}$ with $$e_{1}=\left(\begin{array}{cccc}
0 & -1 & 0 & 0\\
1 & 0 & 0 & 0\\
0 & 0 & 0 & -1\\
0 & 0 & 1 & 0\\
\end{array}\right)
$$
and
$$e_{2}=\left(\begin{array}{cccc}
0 & 0 & 1 & 0\\
0 & 0 & 0 & 1\\
0 & 0 & 0 & 0\\
0 & 0 & 0 & 0\\
\end{array}\right)
$$
satisfies is such that $[D_{\nu}^{*},B]=0$. Consider the Hamiltonian matrix $A$ of the linearised vector field of the R4BP
$$A_{\mu}=\left(\begin{array}{cccc}
0 & 1 & 1 & 0\\
-1 & 0 & 0 & 1\\
a & 0 & 0 & 1\\
0 & b & -1 & 0\\
\end{array}\right)
$$
where $a=U_{x_{1}x_{1}}(0,\mu)$ and $U_{x_{2}x_{2}}(0,\mu)$. We have found a linear symplectic transformation $P$ which brings $A_{\mu_{1}}$ into its normal form $B=P^{-1}A_{\mu_{1}}P$, as we said the matrix
\begin{equation}
\label{versalmat}
B_{\nu}=\left(\
\begin{array}{cccc}
0 & -(1+\nu_{1}) & \nu_{2} & 0\\
(1+\nu_{1}) & 0 & 0 & \nu_{2}\\
-1 & 0 & 0 & -(1+\nu_{1})\\
0 & -1 & 1+\nu_{1} & 0\\
\end{array}\right)
\end{equation}
is a versal deformation of the matrix $B$. Consider now the family of matrices $B_{\mu}=P^{-1}A_{\mu}P$, $B_{\mu}$ is a deformation of $B$ (in $\mu$) so, by the definition
(\ref{def}) there exists a change of parameters $\phi:\Sigma\rightarrow\Lambda$ such that $B_{\mu}=C_{\mu}D_{\nu}C_{\mu}^{-1}$ therefore $B_{\nu}$ and $B_{\mu}$ have the same characteristic polynomials, i.e., the polynomials $P_{B_{\mu}}=\lambda^{4}+(2-a-b)\lambda^{2}+(a+b+ab+1)$ and $P_{B_{\nu}}=\lambda^{4}+2[(1+\nu_{1})^{2}+\nu_{2}]\lambda^{2}+[(1+\nu_{1})^{2}-\nu_{2}]^{2}$
are equal. A straightforward calculation shows: $$\nu_{2}=1/2-1/4(a+b)-1/2(a+b+ab+1)^{1/2}$$ $$\nu_{1}=[1/2-1/4(a+b)+1/2(a+b+ab+1)^{1/2}]^{1/2}-1$$ furthermore $(1+\nu_{1})^{2}=1-1/2(a+b)-\nu_{2}$, the eigenvalues are given by the expression $$\lambda^{2}=-((1+\nu_{1})^{2}-\nu_{2})\pm2\vert1+\nu_{1}\vert\sqrt{\nu_{2}}$$ note that the sign of $\nu_{2}$ determines the behaviour of the eigenvalues i.e., whether they are purely imaginaries or with non-zero real part. In polar coordinates the versal deformation related to the matrix (\ref{versalmat}) looks like $$H_{2,\nu}=\frac{1}{2}\left(R^{2}+\frac{\Theta^{2}}{r^{2}}\right)+\frac{\nu_{2}}{2}r^{2}+(1+\nu_{1})\Theta$$. The parameters $\nu_{1}$ and $\nu_{2}$ are not independent and $\Theta$ is a first integral so we have that the versal deformation depends only of the detuning parameter $\nu_{2}$. Now we are ready to make a study of the truncated system through resonance, the Hamiltonian under consideration will be
\begin{displaymath}
H=\frac{1}{2}\left(R^{2}+\frac{\Theta^{2}}{r^{2}}\right)+\frac{\nu}{2}r^{2}+\Theta+H_{0}^{2}
\end{displaymath}
Because $\Theta$ is a constant we have obtained a system with one degree of freedom, we want to see which terms are the most important near the origin so we will use the scaling of \cite{MeyMcSw} in polar form
\begin{equation}
r\rightarrow\epsilon r
\end{equation}
$$\theta\rightarrow\theta$$ $$R\rightarrow\epsilon^{2}R$$ $$\Theta\rightarrow\epsilon^{3}\Theta$$ $$\nu\rightarrow\epsilon^{2}\nu$$ which is symplectic with multiplier $\epsilon^{-3}$. The Hamiltonian becomes
\begin{equation}
\label{meanterms}
H=\Theta+\epsilon\left(\frac{1}{2}\left(R^{2}+\frac{\Theta^{2}}{r^{2}}\right)+\frac{\nu}{2}r^{2}+h_{1}r^{4}\right)+\mathcal{O}(\epsilon^{2})
\end{equation}
If we ignore the higher order terms, we obtain the form of the truncated Hamiltonian considered in \cite{MeyMcSw} with parameters $\delta=1$ and $\eta=h_{1}>0$. So from the conclusions of that work, we have that the stable and unstable manifolds of the truncated system are globally connected for $\nu<0$ and shrink to the equilibrium as $\nu\rightarrow0^{-}$, actually, when we consider small perturbations of the truncated system or equivalently, when we consider the higher order terms of the full system for $\epsilon$ sufficiently small the stable and unstable manifold may not agree but they will still intersect. This last affirmation is a an application of the Poincar\'e's argument that basically tells us that nondegenerate homoclinic points are stable under small perturbations.
In the references \cite{MeyMcSw} and \cite{McGe} the reader can find detailed discussions on these topics. Therefore, we have satisfied the required conditions under which the ``Blue sky catastrophe" termination occurs: a topologically transverse intersection of the stable and unstable manifolds of the equilibrium point $L_{2}$.

\section{The unstable and stable manifolds}
\begin{figure}[!h]
\begin{center}
\includegraphics[width=4.0in]{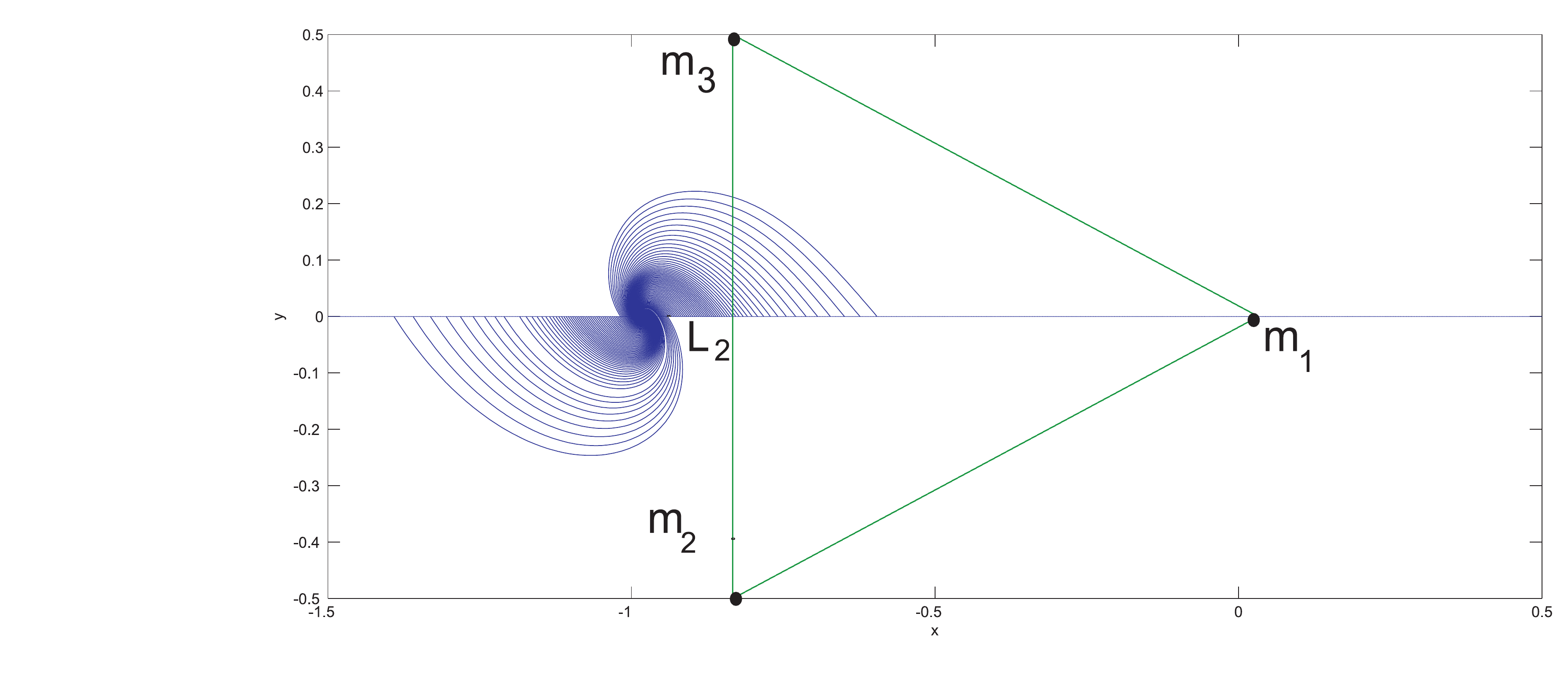}
\end{center}
\caption{The fourth cut of the unstable manifold of $L_{2}$ in the $(x,y)$ plane for $\mu=0.019$.\label{opcuartocorte}}
\end{figure}

\begin{figure}[!h]
\begin{center}
\includegraphics[width=4.5in]{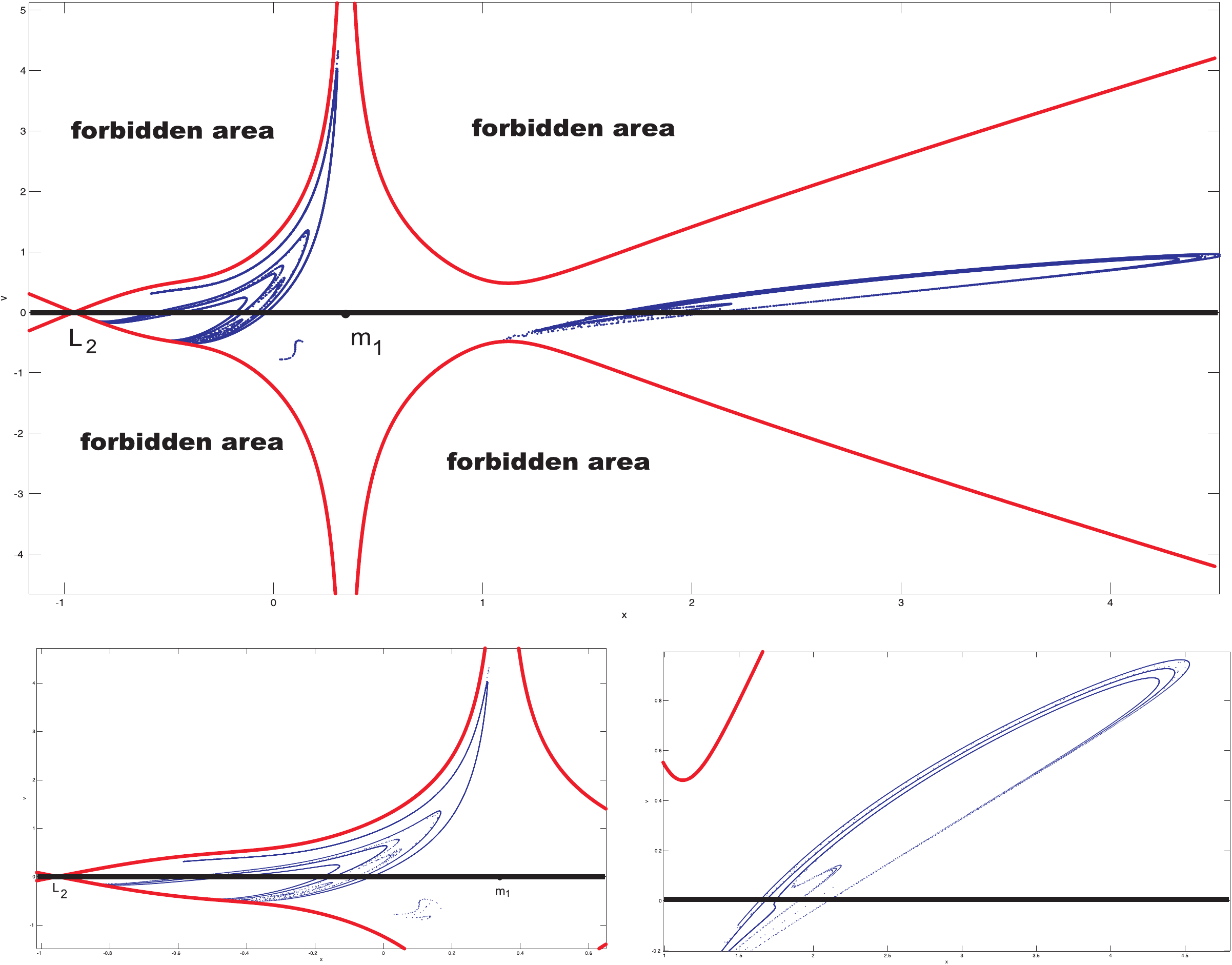}
\end{center}
\caption{First row: the fourth cut of $W^{u}$ for $\mu=0.2$ and $\theta\in[0,\pi]$. Second row: magnifications of the intersections with the $x$-axis}. \label{cuartocortenegativo}
\end{figure}

\begin{figure}[!h]
\begin{center}
\includegraphics[width=4.5in]{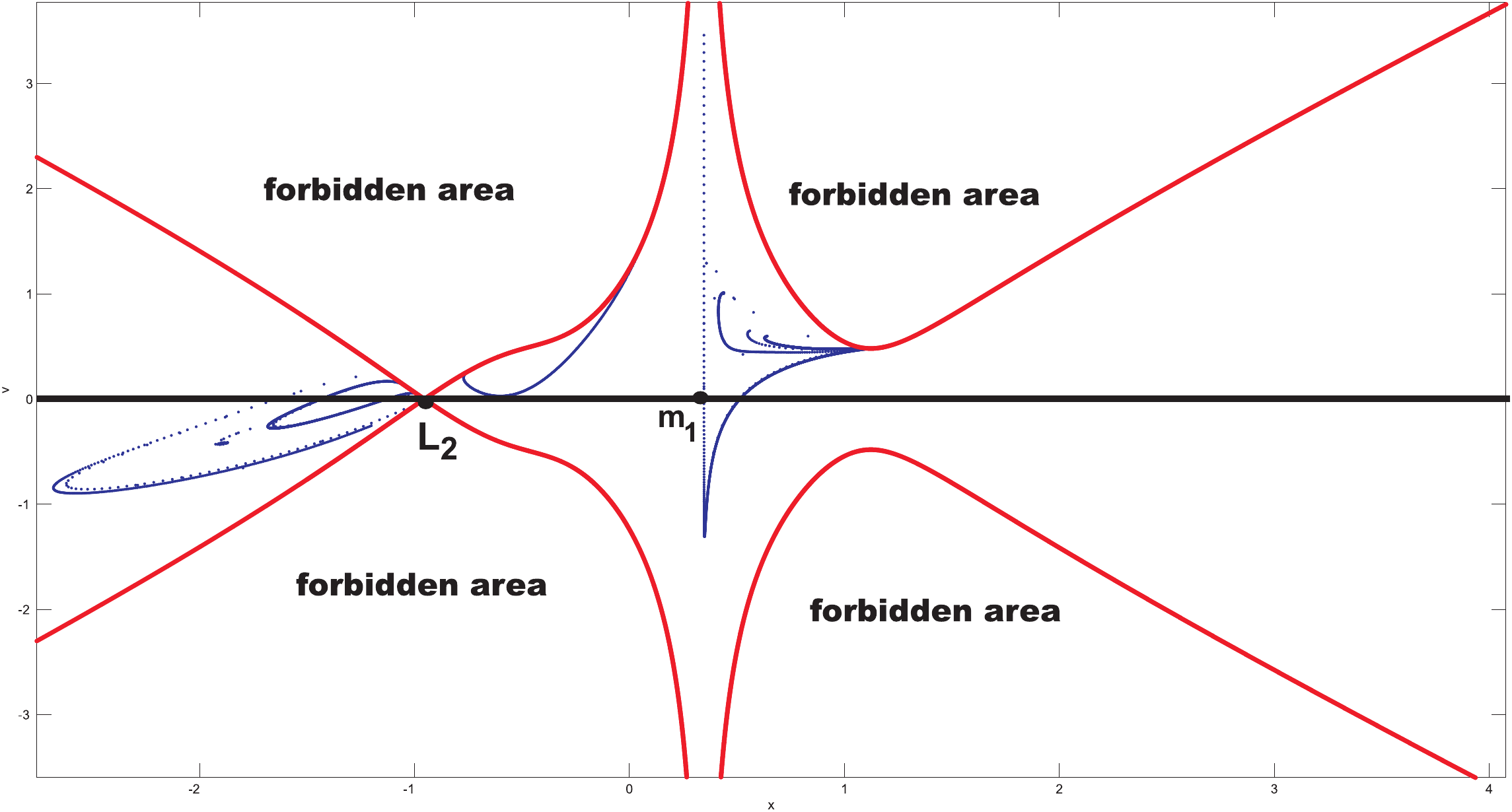}
\end{center}
\caption{The fourth cut of $W^{u}$ for $\mu=0.2$ and $\theta\in[\pi,2\pi]$}.
\label{cuartocortepositivo}
\end{figure}

\begin{figure}
  \centering
   \includegraphics[width=4.5in]{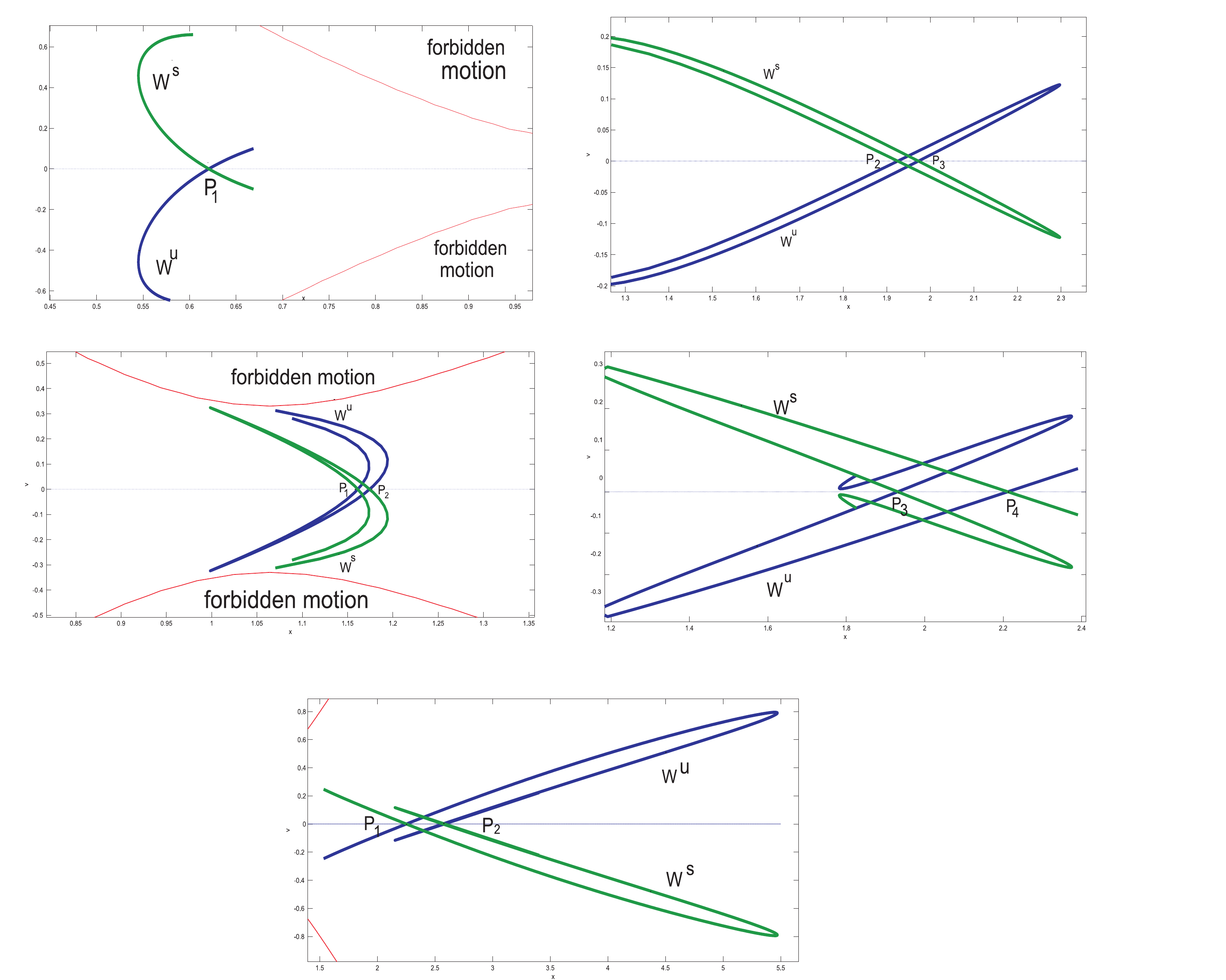}
\caption{First row: intersection points for $\mu=0.019$ at the fifth cut. Second row: intersection points for $\mu=0.1$ at the fourth cut. Third row: intersection points for $\mu=0.2$ at the fourth cut.}\label{cortes}
\end{figure}

In this section we perform a brief numerical study of the unstable manifold $W^{u}$ and the stable manifold $W^{s}$ to show the transverse intersections which were predicted in the former section. We do not attempt to make an exhaustive numerical study of such manifolds so we will restrict ourselves to show some transverse intersections for a few values of the parameter $\mu$ we were able to calculate numerically. We recall that from the linearization around the equilibrium point $L_{2}$ was discussed in the section 3 we saw that for that point there exist a value of the mass parameter $\mu_{b}$ such that for all $\mu\in(\mu_{b},1/3]$ the eigenvalues given by (\ref{charpoly}) are $\lambda_{1,2,3,4}=\pm\alpha\pm\textit{i}\omega$ where $\alpha$ and $\omega$ are real and strictly positive quantities. The general theory states that for this case the solution of the linear part of the equations (\ref{sistemastandar}) is given by
\begin{equation}
\label{ecuacioneslineales} x(t)=\sum_{i=1}^{4}A_{i}e^{\lambda_{i}t}
\end{equation}
$$y(t)=\sum_{i=1}^{4}B_{i}e^{\lambda_{i}t}$$ the coefficients are not independent. In \cite{Sz} it can be found that such coefficients satisfy the relation $(\lambda_{i}^{2}-\Omega_{xx})A_{i}=(2\lambda_{i}+\Omega_{xy})B_{i}$ but the equilibrium point is collinear so $\Omega_{xy}=0$ then $$B_{i}=\left(\frac{\lambda_{i}^{2}-\Omega_{xx}}{2\lambda_{i}}\right)$$ However, as we have seen along this paper, we can not know analytical relations between the eigenvalues and the mass parameter and therefore we do not have a explicit relation between the constants $A_{i}$ ,$B_{i}$ and the mass parameter $\mu$. So for a numerical globalization of the invariant manifolds $W^{u}$ and $W^{s}$ we have decided to take a set of initial conditions totally equivalent to the classical initial conditions shown in \cite{Gomez}. Our vector of initial conditions is given by
\begin{equation}
\label{inicialesnumericas} v_{0}=L_{2}+\epsilon(\cos\theta\bar{v}+\sin\theta\bar{w})
\end{equation}
where $L_{2}=(x_{L_{2}},0,0,0)$, $\epsilon>0$, $\theta\in[0,2\pi]$ and $\bar{v}$ and $\bar{w}$ are in the unstable (stable) space $E^{u}$ $(E^{s})$. As we have already mentioned, the equations of motion have the property that if $x=x(t)$, $y=y(t)$ is a solution, then $x=x(-t)$, $y=-y(-t)$ is also a solution so the unstable and stable manifolds of the equilibrium point are symmetric with respect to the $x-$axis and therefore is enough to calculate the manifold $W^{u}$. In the R3BP and more recently in the restricted four body problem (where the primaries are in a collinear configuration) several authors \cite{Gomez}, \cite{PapaIII} have studied the ``Blue Sky Catastrophe'' termination and the asymptotic orbits by analyzing the transverse intersections of the invariant manifolds. In these cases the equilibrium points, where the BSC termination ocurrs, are off the $x-$axis but in our case the equilibrium point $L_{2}$ is on this axis, so we have to be careful when developing the numerical globalization. As in the R3BP we use the $(x,\dot{x})$ plane to represent the transverse cuts of the invariants manifolds, this can be done because we are using the surface section $y=0$ and because the velocity $\dot{y}$ can be put as function of $x$ and $\dot{x}$ from the relation (\ref{integralprimera}). Using (\ref{sistemastandar}) and (\ref{inicialesnumericas}) we compute and show in the figure (\ref{cortes}) some transverse intersections of $W^{u}$ and $W^{s}$ on the symmetry axis. From the linear approximation we have that the flow is of spiral type in a small neighbourhood of the equilibrium point $L_{2}$ then the orbits on the unstable manifold have several intersections (cuts) with the surface section $y=0$.  Therefore in the "first" cuts we do not expect having a transverse intersection $\dot{x}=0$ of $W^{u}$ in the plane $(x,\dot{x})$, so we have to consider the cuts when the orbits leave an appropriate neighbourhood that depends of the parameter $\epsilon$. In the figure (\ref{opcuartocorte})
we can see that after four cuts the orbits do not have orthogonal intersections ($\dot{x}=0$) with the symmetry axis and therefore the unstable manifold does not intersect the $x-$ axis in the $(x,\dot{x})$ plane, when we continue the integration we observe that the influence of the primaries $m_{2}$ and $m_{3}$ appears on the orbits, if we use a regularization process this is not problem, the interesting phenomena is that the continuity of the $W^{u}$ observed in the first cuts is lost. In the figures (\ref{cuartocortenegativo}) and (\ref{cuartocortepositivo}) we show the unstable manifold for $\mu=0.2$, in order to avoid congestion of the picture we present the manifold for $\theta\in [0,\pi]$ and $\theta\in [\pi,2\pi]$.

\begin{figure}[!h]
\begin{center}
\includegraphics[width=4.5in]{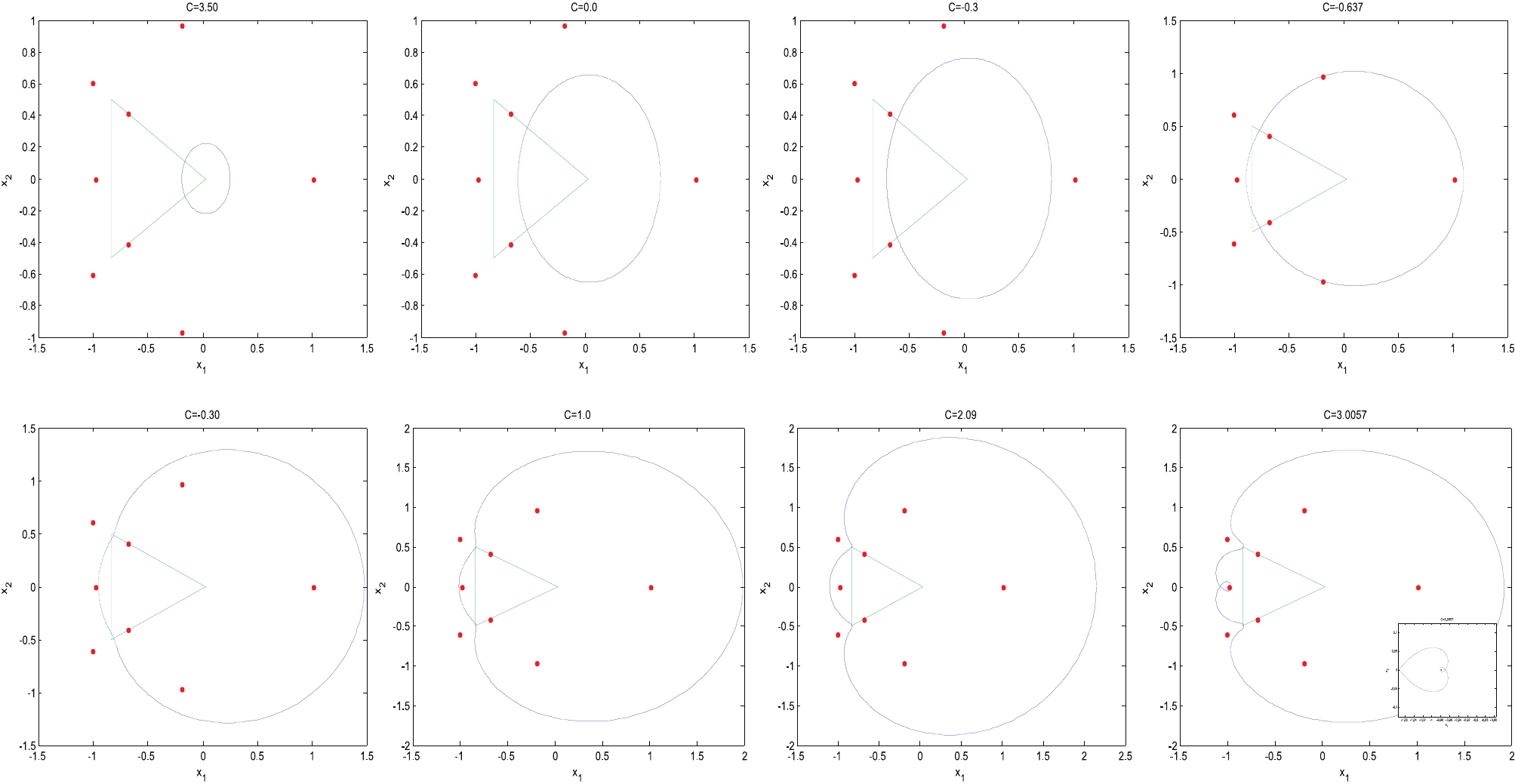}
\end{center}
\caption{The blue sky catastrophe termination for the family $f$.\label{fasesf}}
\end{figure}

However we can find transverse intersections of $W^{u}$ on the symmetry axis. For the value of the parameter $\mu=0.019$ we found the transversal intersection $P_{2}\approx1.925$ corresponding to the homoclinic orbit that is the termination of the family $f$ \cite{Burgos}, see figure (\ref{fasesf}). If the value of $\mu$ is increased, more distant transverse intersections of the three primaries appear, see figure (\ref{cortes}).
\newpage
\textbf{Acknowledgements} The authors thank Jaume Llibre for his helpful comments and suggestions. Author Burgos--Garc\'ia has been supported by a CONACYT fellowship of doctoral studies.

\end{document}